\documentclass[11pt]{article}
\usepackage{amssymb,amsfonts,amsmath,latexsym,epsf,tikz,url}

\newtheorem{theorem}{Theorem}[section]
\newtheorem{proposition}[theorem]{Proposition}
\newtheorem{observation}[theorem]{Observation}

\newtheorem{corollary}[theorem]{Corollary}
\newtheorem{lemma}[theorem]{Lemma}

\newtheorem{definition}[theorem]{Definition}

\newcommand{\proof}{\noindent{\bf Proof. }}
\newcommand{\qed}{\hfill $\square$\medskip}

\begin{document}

\title{On the edge chromatic vertex stability number  of graphs}

\author{
Saeid Alikhani$^{}$\footnote{Corresponding author}
\and
Mohammad R. Piri
}

\date{\today}

\maketitle

\begin{center}
Department of Mathematics, Yazd University, 89195-741, Yazd, Iran\\
{\tt alikhani@yazd.ac.ir, piri4299@gmail.com}
\end{center}


\begin{abstract}
For an arbitrary invariant $\rho (G)$ of a graph $G$, the $\rho-$vertex stability number $vs_{\rho}(G)$ is the minimum number of vertices of $G$ whose removal results in a graph $H\subseteq G$ with $\rho (H)\neq \rho (G)$ or with $E(H)=\varnothing$. In this paper,
first  we give some general lower and upper bounds for the $\rho$-vertex  stability number, and then study the edge chromatic  stability number of graphs, $vs_{\chi^{\prime}}(G)$,  where $\chi^{\prime}=\chi^{\prime}(G)$ is edge chromatic number (chromatic index) of $G$. We prove some general results for this parameter  and determine $vs_{\chi^{\prime}}(G)$  for specific classes of graphs.
\end{abstract}


\medskip
\noindent{\bf AMS Subj.\ Class.}: 05C15, 05C25. 

\section{Introduction and definitions}
Let $G = (V,E)$ be a simple graph with $n$ vertices. Throughout this paper we consider only simple graphs.  A graph is
empty if $E(G) = \varnothing$. We follow the standard graph notation and definitions \cite{DMGT}.
A (graph) invariant $\rho(G)$ is a function $\rho : {\cal I}\longrightarrow \mathbb{R}^+\cup \{\infty \}$, where
${\cal I}$ is the class of finite simple graphs. An invariant $\rho(G)$ is integer valued if its
image set consists of non-negative integers, that is, $\rho({\cal I})\subseteq \mathbb{N}_0$.
An invariant $\rho(G)$ is monotone increasing if $H \subseteq G$ implies $\rho(H)\leq \rho(G)$,
and monotone decreasing if $H\subseteq G$ implies $\rho(H)\geq \rho(G)$; $\rho(G)$ is monotone if it
is monotone increasing or monotone decreasing.
If $H_1$ and $H_2$ are disjoint graphs, then an invariant is called additive if $\rho (H_1 \cup H_2)=\rho(H_1) 
+ \rho(H_2)$ and maxing if $\rho (H_1 \cup H_2)= max \{\rho (H_1), \rho (H_2)\}$.

Kemnitz and Marangio in \cite{DMGT}, for an arbitrary invariant of $G$, introduced   the $\rho-$edge  stability number $es_{\rho}(G)$ of $G$ as   the minimum number of edges  of $G$ whose
removal results in a graph $H\subseteq G$ with $\rho (H) \neq \rho (G)$ or with $E(H)=\emptyset$. Then they 
gave some general lower and upper bounds for the $\rho$-edge stability number. Also they studied  the $\chi'$-edge stability number of graphs, where $\chi'=\chi'(G)$ is the chromatic index of $G$. They proved some general results for the so-called chromatic edge stability index $es_{\chi'}(G)$ and determine $es_{\chi'}(G)$ exactly for specific classes of graphs. For recent results on chromatic edge stability number and stabilizing on the distinguishing number of a graph, see \cite{COMM,Henning,Nazanin}. 
  
  \medskip

Motivated  by these papers and especially  \cite{DMGT}, we state the following definitions:

\begin{definition}
For an arbitrary invariant $\rho$ of $G$,  the $\rho-$vertex stability number $vs_{\rho}(G)$ of $G$ is  the minimum number of vertices of $G$, whose
removal results in a graph $H\subseteq G$ with $\rho (H) \neq \rho (G)$ or with $E(H)=\emptyset$.
\end{definition}
 
\begin{definition}
The chromatic vertex stability number $vs_{\chi^{\prime}}(G)$ of a non-empty graph $G$, is the minimum number of vertices of $G$ whose removal results in graph $H\subseteq G$ with $\chi^{\prime}(H)\neq \chi^{\prime} (G)$. If $G$ is empty, then $vs_{\chi^{\prime}}(G)=0$.
\end{definition}

In this paper, we first consider the general case and give bounds for arbitrary $\rho$-vertex stability numbers of graphs. In Section 3, we  study  $\chi'$-vertex stability number of graphs, and compute the $vs_{\chi'}(G)$ for specific graphs.  In the last section,  we study the $\chi'$-vertex stability number of join and corona product of two graphs.

\section{On the $\rho$-vertex stability number}
In this section, we obtain some bounds for arbitrary  $\rho$-vertex stability number of graphs. We follow \cite{DMGT} and  similar to $es_{\rho}(G)$, we obtain results for $vs_{\rho}(G)$.

\begin{theorem}\label{aaa}
Let $\rho(G)$ be maxing and monotone increasing, let $G=H_1 \cup H_2 \cup ...  \cup H_k$ be a graph whose subgraphs $H_1, ... , H_k$ and the integer $s\geq 1$ are defined such that $\rho (H_i)=\rho (G)$ if and only if $1\leq i \leq s$. Then $vs_{\rho}(G)=\sum_{i=1}^s vs_{\rho}(H_i)$.
\end{theorem}
\proof
Let $V^\prime =V_1^\prime \cup ... \cup V^\prime _s$ with $V^\prime _i \subseteq V(H)$, $|V^\prime _i|= vs_\rho (H_i)$ and $\rho (H_i - V^\prime _i ) \neq \rho (H_i)$  for $i=1, ... ,s$. Since the invariant is maxim, $max \{\rho (H_i - V^\prime _i ): 1\leq i \leq s\} \cup \{ \rho(H_i): s+1\leq i \leq k\} \neq \rho (G)$ which implies $vs_\rho (G) \leq \mid V^\prime \mid =\sum_{i=1}^s vs_\rho (H_i)$. If a vertex set $V^{\prime\prime}$ with less than $\mid V^\prime \mid$ vertex  is removed from $G$, then there is a subgraph $H_j$, $1\leq j\leq s$, from which less than $vs_\rho (H_j)$ vertices are removed, which implies $\rho (H_j - V^{\prime\prime} )=\rho(H_j)$ and thus, since the invariant is maxing and monotone increasing, $\rho (G-V^{\prime\prime})=\rho(H_j)=\rho(G)$. Therefore, $vs_\rho (G)=|V^\prime| =\sum_{i=1}^s vs_\rho (H_i)$.\qed

\begin{theorem}\label{aa}
Let $\rho(G)$ be monotone and let $G$ be a nonempty graph with $\rho(G)=k$. If $G$ contains $s$ nonempty subgraphs $G_1, ... , G_s$ with $\rho(G_1)= ... = \rho(G_s)= k$ such that $a\geq 0$ is the number of vertices that occur in at least two of these subgraphs and $q\geq 1$ is the maximum number of these subgraphs with a common vertex, then both $vs_\rho(G)\geq \frac{1}{q} \sum_{i=1}^s vs_\rho (G_i) \geq s/q$ and $vs_\rho(G)\geq \sum_{i=1} ^s vs_\rho (G_i)-a(q-1)$ hold. 
\end{theorem}
\proof
Let $\rho(G)$ be monotone increasing. Let $V^\prime$ be a set of vertices of $G$ with $|V^\prime|=vs_\rho(G)$ such that $\rho(G-V^\prime) <k$. Then the set $V^\prime$ must contain at least $vs_\rho(G_i)$ vertices of $G_i$, $1\leq i \leq s$, since otherwise $k> \rho (G-V^\prime )\geq \rho (G_i -V^\prime \cap V(G_j))=k$ for some $j$, $1\leq j \leq s$, a contradiction. Therefore, $b=\sum_{i=1} ^s \mid V^\prime \cap V(G_i)\mid \geqslant \sum_{i=1} ^s vs_\rho(G_i) \geqslant s$.\\
On the other hand, at most $a^\prime =min\{a,|V^\prime|\}$ vertices of $V^\prime$ are counted at most $q$ times in $b$, every other vertices of $V^\prime$ is counted at most once, so $b\leqslant a^\prime .q + (|V^\prime|-a^\prime ).1= |V^\prime|+ a^\prime (q-1)$. Since $a^\prime \leqslant |V^\prime|$, $b\leqslant q |V^\prime|$  so, $vs_\rho (G) =| V^\prime|\geqslant b/q \geqslant 1/q \sum_{i=1} ^s vs_\rho (G_i) \geqslant s/q$. Also, $a^\prime \leqslant a$ implies $vs_\rho (G) = \mid V^\prime \mid \geqslant b-a (q-1) \geqslant \sum_{i=1} ^s vs_\rho (G_i) - a(q-1)$. 
The proof for monotone decreasing $\rho (G)$ runs analogously.\qed
\begin{corollary}\label{bb}
Let $\rho (G)$ be monotone and let $G$ be a nonempty graph with $\rho (G)=k$. If $G$ contains $s$ nonempty subgraphs $G_1, ... , G_s$ with $\rho (G_1)= ... = \rho (G_s)= k$ and pairwise disjoint vertices  sets, then $vs_\rho (G) \geqslant \sum_{i=1} ^s vs_\rho (G_i) \geqslant s$.
\end{corollary}
\proof
Each vertex of $G$ is contained in at most $q=1$ of the given subgraphs since they are are pairwise vertex disjoint. The results follows from Theorem \ref{aa}. \qed
\begin{corollary}\label{bbb}
Let $\rho (G)$ be monotone. If $H\subseteq G$ and $\rho (H)=\rho (G)$, then $vs_\rho (H) \leqslant vs_\rho (G)$.
\end{corollary}
\proof
Corollary \ref{bb} with $s=1$ implies the result.\qed
\begin{proposition}
If $G$ is a graph without isolated vertex, then $vs_{\delta}(G)=1$.
\end{proposition}
\proof
It suffices to remove one vertex adjacent with a vertex of degree $\delta (G)$ in order to decrease the minimum degree, hence $vs_\delta (G)=1$.\qed

\medskip 

Here we state the following definition, to obtain more results:   
\begin{definition}
Let $V_\Delta$ be the set of vertices of $G$ of degree $\Delta (G)$. A dominating set for $V_\Delta$ is a subset $\Gamma$ of $V$ such that $N_G (\Gamma )=V_\Delta$. The domination number $\gamma (V_\Delta )$ is the minimum number of vertices in a dominating set for $V_\Delta$.
\end{definition}
\begin{proposition}\label{cc}
For a connected graph $G$, $vs_\Delta (G)=\gamma (V_\Delta )$.
\end{proposition}
\proof
If $G$ is not empty, then $\Delta (G)\geq 1$. Let $\Gamma$ be minimum dominating  set of vertices of $V_\Delta$ that $\mid \Gamma \mid =\gamma (V_\Delta )$. So $\Delta (G-\Gamma )< \Delta (G)$ and since $\Gamma$ is minimum, we have  $vs_\Delta (G)=\gamma (V_\Delta)$.\qed

\medskip

Let $\omega(G)$ be the number of components of a graph $G$ and $\kappa(G)$ the vertex connectivity of $G$, that is, the minimum number of vertices whose removal gives a disconnected graph or the singleton $K_1$. By the definition it follows that if $G$ is connected, then $vs_{\omega}(G)=\kappa(G)$.
\begin{proposition}
Let $G$ be a graph with $\omega(G)$ components $H_1, ... ,H_{\omega(G)}$. Then $vs_{\omega} (G)=min \{ \kappa (H_i) : 1\leq i \leq \omega (G), H_i \ncong K_1\}$ if $G$ is not empty.
\end{proposition}
\proof
The number of components is additive and can be increased by vertex deletions for nonempty graphs. Let $H_1, ... , H_{\omega(G)}$ be the nonempty components of $G$. Then $vs_{\omega} (G)= min \{ vs_{\omega} (H_i) : 1\leq i \leq k(G) \} = min \{ \kappa (H_i) : 1\leq i \leq \omega(G) \}$. \qed

\section{General results for the chromatic vertex stability index}
A function $c : E(G) \longrightarrow \{1, . . . , k\}$ such that
$c(e_1) \neq c(e_2)$ for any two adjacent edges $e_1$ and $e_2$ is called a $k$-edge coloring
of $G$, and $G$ is called $k$-edge colorable. The minimum $k$ for which $G$ is $k$-edge
colorable is the chromatic index $\chi^\prime (G)$ of $G$. By Vizing's Theorem, the chromatic
index can only attain one of two values, $\Delta (G) \leq \chi^\prime (G) \leq \Delta (G)+1$. Graphs with
$\chi^\prime (G)=\Delta (G)$ are called class $1$ graphs and graphs with $\chi^\prime (G)=\Delta (G)+1$ are
called class $2$ graphs. As in \cite{DMGT}, we consider  the invariant $class(G) = \chi^\prime (G)-\Delta(G)+1 \in \{1, 2\}$.
A graph $G$ is called overfull if its order $n$ is odd and if it contains more than
$\Delta (G) (n-1)/2$ edges. Obviously, an overfull graph must be a class $2$ graph.
 Note that $\chi^\prime (G)$ is an invariant which is monotone increasing, integer valued,
and maxing \cite{DMGT}.

\medskip
In this section, we consider the $\chi^\prime$-vertex  stability number, $vs_{\chi^\prime}(G)$, which we
also call chromatic vertex stability index of $G$. Using Theorem \ref{aaa} we can compute
$vs_{\chi^\prime}(G)$ by the chromatic vertex stability indices  of its components. Let $G = H_1 \cup ... \cup H_k$ such that $\chi^\prime(G)=\chi^\prime(H_i)$ if and only if $1 \leq i \leq s$ for $s \leq k(G)$. Then
$vs_{\chi^\prime}(G)=\sum_{i=1} ^s vs_{\chi^\prime}(H_i)$. Therefore, we can assume without loss of generality in
the following that $G$ is connected.
In \cite{DMGT} proved that $es_{\chi^\prime}(G) \leq \lfloor \vert E(G)\vert / \chi^\prime (G) \rfloor \leq \alpha^\prime (G)$, and since $vs_{\chi'}(G) \leq es_{\chi'}(G)$, so  $vs_{\chi^\prime}(G) \leq \lfloor \vert E(G)\vert / \chi^\prime (G)\rfloor$.

\begin{lemma}\label{dd}
If $G$ is a class $1$ graph, then $vs_{\chi^\prime}(G)\geq vs_\Delta (G)$.
\end{lemma}
\proof
If $G$ is nonempty, then there is a set $V^\prime$ of vertices of $G$ such that $\mid V^\prime \mid =vs_{\chi^\prime}(G)$ and $\Delta (G-V^\prime ) \leq \chi^\prime (G-V^\prime ) < \chi^\prime (G) = \Delta (G)$. It follows that $\Delta (G-V^\prime ) < \Delta (G)$ which implies $vs_{\chi^\prime}(G) =\mid V^\prime \mid \geq vs_\Delta (G)$.\qed
\begin{proposition}\label{ee}
If $G$ is a class $1$ graph and there is a vertex set $V^\prime$ such that $\mid V^ \prime \mid = vs_\Delta(G)$, $\Delta (G-V^\prime ) < \Delta (G)$ and $(G-V^\prime )$ is in the class $1$, then $vs_{\chi^\prime}(G) = vs_\Delta (G)$.
\end{proposition}
\proof
By Lemma \ref{dd}, $vs_{\chi^\prime}(G)\geq vs_\Delta (G)$ and from the properties of the given set $V^\prime$, since $\chi^\prime (G-V^\prime ) = \Delta (G-V^\prime ) < \Delta (G) = \chi^\prime (G)$ which implies $vs_\Delta (G) \geq vs_{\chi^\prime}(G)$.\qed

\medskip  
Proposition \ref{ee} and Lemma \ref{dd} imply that if $G$ is in the class $1$ but $G-V^\prime$ is in the class $2$ for all sets with $\mid V^\prime \mid = vs_\Delta (G)$ and $\Delta (G-V^\prime ) < \Delta (G)$, then $vs_{\chi^\prime} (G) >  vs_\Delta (G)$. An example for such graphs is the complete graph $K_4$.
\begin{theorem}
If $G$ is a class $2$ graph, then $vs_{\chi^\prime}(G) = min \{ vs_\Delta (G), vs_{class}(G) \}$.
\end{theorem}
\proof
Since $G$ is in the class $2$, the graph $G$ is not empty and invariants $\Delta (G)$, $class (G)=2$ and $\chi^\prime (G)=\Delta (G)+1$ can be reduced by vertex removal. By removing $vs_\Delta (G)$ vertex $V^\prime$ such that $\Delta (G-V^\prime ) < \Delta (G)$ we obtain $\chi^\prime (G-V^\prime ) \leq \Delta (G-V^\prime ) +1 < \Delta (G) +1=\chi^\prime (G)$ which implies $\mid V^\prime \mid \geq vs_{\chi^\prime} (G)$. By removing $vs_{class}(G)$ vertices $V^{\prime\prime}$ such that $class (G-V^{\prime\prime})=1$ we obtain $\chi^\prime (G-V^{\prime\prime})=\Delta (G-V^{\prime\prime}) \leq \Delta (G) < \Delta (G)+1=\chi^\prime(G)$ which implies $\mid V^{\prime\prime} \mid \geq vs_{\chi^\prime}(G)$. It follows that $min \{ vs_\Delta (G), vs_{class}(G) \} \geq vs_{\chi^\prime}(G)$.\\
Consider now a set of vertices $(G-V^{\prime\prime\prime})$ such that $\chi^\prime (G-V^{\prime\prime\prime}) < \chi^\prime (G)=\Delta (G)+1$, that is, $\chi^\prime (G-V^{\prime\prime\prime}) \leq \Delta (G)$. Then $G-V^{\prime\prime\prime}$ cannot both be in class $2$ and have the same maximum degree as $G$ since this would imply $\chi^\prime (G-V^{\prime\prime\prime})=\Delta (G)+1$. Therefore, $|V'''| \geq vs_\Delta (G)$ or $\mid V^{\prime\prime\prime} \mid \geq vs_{class} (G)$ which implies $vs_{\chi^\prime}(G) \geq min \{ vs_\Delta (G), vs_{class}(G) \}$.\qed
\begin{theorem}
If $G$ is bipartite, then $vs_{\chi^\prime}(G)=vs_\Delta (G)$.
\end{theorem}
\proof
The result follows from Proposition \ref{ee} since every subgraph $G-V^\prime$ of $G$ is bipartite and thus in class $1$.\qed
\begin{corollary}
 $ 	vs_{\chi^\prime} (K_{m,n})=vs_\Delta (K_{m,n})=\left\{
  	\begin{array}{ll}
  	{\displaystyle
  		2}&
  	\quad\mbox{if $m=n$, }\\[15pt]
  	{\displaystyle
  		1} &
  	\quad\mbox{if $m\neq n$. }
  	\end{array}
  	\right.	
  	$
\end{corollary}
\proof
We know that $\chi^\prime (K_{m,n})=max\{m,n\}$. If $m\neq n$ and $m>n$, then by removing a vertex of degree $m$, we have $K_{m-1,n}$ that $\chi^\prime (K_{m-1,n})=m-1$. So in this case, $vs_{\chi^\prime}(K_{m,n})=1$.\\
Now if $m=n$, by removing  a vertex, we have $K_{m,m-1}$, that $\chi^\prime (K_{m,m-1})=m$. So $vs_{\chi^\prime}(K_{m,m})>1$. By the previous part, $vs_{\chi^\prime}(K_{m,m-1})=1$. Therefore $vs_{\chi^\prime}(K_{m,m})=vs_\Delta (K_{m,m})=2$. \qed
\begin{observation}\label{ff}
If $\chi^\prime (G)=2$, then $vs_{\chi^\prime}(G)=vs_\Delta (G)$.
\end{observation}
\begin{theorem}\label{gg}
For $n\geq 3$, $vs_{\chi^\prime}(P_n)=vs_\Delta (P_n)=\gamma (P_{n-2})=\lceil \frac{n-2}{3} \rceil$.
\end{theorem}
\proof
By Observation \ref{ff}, $vs_{\chi^\prime}(P_n)=vs_\Delta (P_n)$. On the other hand, by Proposition \ref{cc}, $vs_\Delta (P_n)=\gamma (V_\Delta )=\gamma (P_{n-2})=\lceil \frac{n-2}{3} \rceil$.\qed

Let $t^*(G)$ be the minimum number of edges in a color class of the graph
$G$ where the minimum is taken over all edge colorings of $G$ with $\chi^\prime(G)$ colors.
If $G$ is nonempty, then removing a vertex of any edges in one color class from $G$ reduces the chromatic index, 
thus $vs_{\chi^\prime}(G)\leq t^*(G)$ follows (\cite{DMGT}).

\begin{theorem}\label{ddd}
For $n\geq 3$, $vs_{\chi^\prime}(C_n)=1$ if $n$ is odd, and $vs_{\chi^\prime}(C_n)=\lceil \frac{n}{3} \rceil$ if $n$ is even.
\end{theorem}
\proof
Since $C_n$ is nonempty, $vs_{\chi^\prime}(C_n) \geq 1$. If $n$ is odd, then $t^* (C_n)=1$. So $vs_{\chi^\prime}(C_n)=t^* (C_n)=1$. If $n$ is  even, then $\chi^\prime (C_n)=2$. By Observation \ref{ff} and Proposition \ref{cc},
$vs_{\chi^\prime}(C_n)=vs_\Delta (C_n)=\gamma (C_n)=\lceil \frac{n}{3} \rceil$. \qed
\begin{theorem}
	\begin{enumerate}
		\item[(i)] 
		There exists a graph $G$ such that $\mid vs_{\chi^\prime}(G)-t^* (G) \mid$ can be arbitrarily large.
		\item[(ii)] 
		There exists a graph $G$ such that $|vs_{\chi^\prime}(G)-vs_{\Delta}(G)|$ can be arbitrarily large.
		\end{enumerate}
		\end{theorem}
\proof
	\begin{enumerate}
		\item[(i)] 
Consider the graph $P_n$. Note that $t^*(P_n)=\lfloor \frac{n}{2} \rfloor$ and by Theorem \ref{gg}, $vs_{\chi^\prime}(P_n)=\lceil \frac{n-2}{3} \rceil$.
	\item[(ii)] 
	Consider the graph $C_n$ with $n$ odd. We have $vs_{\chi^\prime}(C_n)=1$ and by Theorem \ref{ddd}, $vs_{\Delta}(C_n)=\lceil \frac{n}{3} \rceil$.\qed
		\end{enumerate}

\begin{theorem}\label{nordhaus}
Suppose that $G$ and $\overline{G}$ are nonempty graphs with sizes $m$ and $\overline{m}$, respectively. Then
\begin{equation*}
2 \leq vs_{\chi^\prime}(G) + vs_{\chi^\prime}(\overline{G}) \leq \lfloor \frac{m}{\chi^\prime (G)} \rfloor + \lfloor \frac{\overline{m}}{\chi^\prime (\overline{G})} \rfloor.
\end{equation*}
\end{theorem}
\proof
Since $G$ and $\overline{G}$ are nonempty, so $vs_{\chi^\prime}(G) \geq 1$ and $vs_{\chi^\prime}(\overline{G}) \geq 1$. Therefore $2 \leq vs_{\chi^\prime}(G) + vs_{\chi^\prime}(\overline{G})$. Let $C_1, ... , C_t$ where $t=\chi^\prime(G)$, be color classes of edge coloring of $G$ and $t^* (G)$ be cardinality  of minimum color class. Because $t^* (G) \leq |C_i|$ for $1 \leq i \leq t$, hence $m=|C_1|+|C_2|+...+
|C_t| \geq t^* (G) \chi^\prime (G)$. Therefore $t^* (G) \leq \lfloor \frac{m}{\chi^\prime (G)} \rfloor$. On the other hand, $vs_{\chi^\prime}(G) \leq t^* (G)$. So $vs_{\chi^\prime}(G) \leq \lfloor \frac{m}{\chi^\prime (G)} \rfloor$. With similar argument, $vs_{\chi^\prime}(\overline{G}) \leq \lfloor \frac{\overline{m}}{\chi^\prime (\overline{G})} \rfloor$. Therefore $vs_{\chi^\prime}(G) + vs_{\chi^\prime}(\overline{G}) \leq \lfloor \frac{m}{\chi^\prime (G)} \rfloor + \lfloor \frac{\overline{m}}{\chi^\prime (\overline{G})} \rfloor$.\qed

\medskip 
The cycle graph $C_4$ is an example of graph for which the equality of upper bound in Theorem \ref{nordhaus} hold.

\begin{theorem}\label{gg}
If $G$ is a graph of order $n$ and $\Delta (G)=n-1$, then $vs_{\chi^\prime}(G)=1$ or $vs_{\chi^\prime}(G)=2$.
\end{theorem}
\proof
First we state and prove two claims: 
\begin{enumerate}
\item[Calim 1)] If $G$ is in the class $2$, then $vs_{\chi^\prime}(G)=1$.
\\ For the proof, suppose that $deg(v)=n-1$ and so $\Delta (G-v) \leq n-2$. So $\chi^\prime (G-v) \leq n-1 < n = \chi^\prime (G)$ and therefore $vs_{\chi^\prime}(G)=1$.
\item[Claim 2)] If we have only one vertex of degree $n-1$, then $vs_{\chi^\prime}(G)=1$.\\
For the proof, suppose that $deg(v)=n-1$ and $deg(u) < n-1$ for every $u \neq v$. Since $\chi^\prime (G)=n$ or $\chi^\prime (G)=n-1$. With removing $v$, we have $\Delta (G-v)=n-1$. So $\chi^\prime (G-v) \leq \Delta (G-v)+1 < n-2+1=n-1 \leq \chi^\prime (G)$. Therefore $vs_{\chi^\prime}(G)=1$
\end{enumerate}
Now suppose $deg(v)=n-1$. So $\Delta (G-v) \leq n-2$. If $G-v$ is in the class $1$, then $vs_{\chi^\prime}(G)=1$. If $G-v$ is in the class $2$ and $\Delta (G-v)=n-2$, then by Claim $1$, $vs_{\chi^\prime}(G-v)=1$. Thus $vs_{\chi^\prime}(G)=2$.\qed

\medskip
A wheel graph $W_n$ with $n\geq 3$ is the join of a cycle $C_n$, with  consecutive vertices $v_1, v_2, ... , v_n$ and a single vertex $w$. Wheels are class $1$ graphs. 
By Theorem \ref{gg} we have the following result: 
\begin{theorem}\label{hh}
\begin{enumerate}
\item[i)]
$vs_{\chi^\prime} (K_n)=\left\{
  	\begin{array}{ll}
  	{\displaystyle
  		1}
  	\quad\mbox{if $n$ is odd, }\\[15pt]
  	{\displaystyle
  		2} 
  	\quad\mbox{if $n$ is even. }
  	\end{array}
  	\right.$	
  
\item[ii)]
 $vs_{\chi^\prime}(W_3)=2$ and for $n\geq 4$, $vs_{\chi^\prime}(W_n)=1$.
\end{enumerate}
\end{theorem}

\medskip 
Let $p$ and $a_1\leq a_2\leq ... \leq a_p$ be positive integers. The complete multipartite graph $K_{a_1, a_2, ... ,a_p}$ is defined as follows. It has $n=a_1+a_2+...+a_p$ vertices, partitioned into parts $A_1, A_2, ... , A_p$ where each $A_i$ has cardinality $a_i$. It two vertices are in the same part, they are not adjacent, while if they are in different parts, they are joined by exactly one edge.

\begin{theorem}{\rm\cite{hr}}
If $G$ is a complete multipartite graph, then
\[
 	\chi^\prime (G)=\left\{
  	\begin{array}{ll}
  	{\displaystyle
  		\Delta (G)}&
  	\quad\mbox{if $G$ is not overfull, }\\[15pt]
  	{\displaystyle
  		\Delta (G)+1} &
  	\quad\mbox{if $G$ is overfull. }
  	\end{array}
  	\right.	
  	\]
\end{theorem}
By  Theorem \ref{gg} we have the following corollary: 
\begin{corollary}
Let $G$ be complete multipartite $K_{1,a_1, ... ,a_p}$ that $2\leq a_2\leq ... \leq a_p$, then $vs_{\chi^\prime}(G)=1$.
\end{corollary}

\begin{theorem}
Let $G$ be complete multipartite  graph $K_{a_1, a_2, ... , a_p}$ that $p\geq 3$. Then $1 \leq vs_{\chi^\prime}(G)\leq 3$.
\end{theorem}
\proof
Let $a_1=a_2= ... =a_p$. If $a_1=1$, then $K_{a_1, a_2, ... , a_p}=K_n$ and  by Theorem \ref{ff}, $vs_{\chi^\prime}(G)=1$ or $vs_{\chi^\prime}(G)=2$. Now suppose $a_1=a_2= ... =a_p \geq 2$. We remove the arbitrary vertex $v$. So $\Delta (G-v)=\Delta (G)$. This concludes that $vs_{\chi^\prime}(G)\geq 1$. Now with removing  a vertex $u$ of maximum degree in $G-v$, $\Delta (G-v-u)=\Delta (G)-1 < \Delta (G)$. We consider two cases.
\begin{enumerate}
\item[Case 1)] If $G$ is not overfull, then  $\chi^\prime (G)=\Delta (G)$. Therefore $vs_{\chi^\prime}(G)=2$, if $\chi^\prime (G-v-u)=\Delta(G)-1$ and $vs_{\chi^\prime}(G)=3$, if $\chi^\prime (G-v-u)=\Delta (G)$.
\item[Case 2)] If  $G$ is overfull, then $\chi^\prime (G)=\Delta (G)+1$. If $\chi^\prime (G-v)=\Delta (G)$, then $vs_{\chi^\prime}(G)=1$. Otherwise $vs_{\chi^\prime}(G)=2$.
\end{enumerate}
Consider $G=K_{a_1, ... ,a_p}$ that $a_i\neq a_j$ for some $1\leq i,j \leq p$. We have the following cases:
\begin{enumerate}
\item[Case a)] Let $\chi^\prime (G)=\Delta (G)$. By removing  a vertex $v$ of part $A_p$, we have $\Delta (G-v)=\Delta (G)-1$. In this case, if $\chi^\prime (G-v)=\Delta (G)-1$, then $vs_{\chi^\prime}(G)=1$ and if $\chi^\prime (G-v)=\Delta (G)$, then by removing a vertex  $u$ from greatest part of $G-v$, we have $\Delta (G-v-u)=\Delta (G)-2$. Thus $vs_{\chi^\prime}(G)=2$.
\item[Case b)] Let $\chi^\prime (G)=\Delta (G)+1$. By removing  a vertex $v$ of part $A_p$, we have $\Delta (G-v)=\Delta(G)-1$. So $\chi^\prime (G-v)=\Delta (G)$ or $\chi^\prime (G-v)=\Delta (G)-1$. It follows that $vs_{\chi^\prime}(G-v)<vs_{\chi^\prime}(G)$. Therefore $vs_{\chi^\prime}(G)=1$.
\end{enumerate}\qed
\begin{theorem}
Let $G$ be a graph which consists of a complete graph $K_n$, $n\geq 2$, and an additional vertex $w$ connected to $d=d(w)$ vertices of $K_n$. Then we have the following.
\begin{enumerate}
\item[(i)]
$vs_\Delta (G)=1$ and if $n$ is even, then for $d=0$, $vs_{\chi^\prime}(G)=2$, and for $1\leq d \leq n$,   $vs_{\chi^\prime}(G)=1$.
\item[(ii)]
If $n\geq 2$ is odd and $0 \leq d \leq n$, then $vs_{\chi^\prime}(G)=2$. 
\end{enumerate}
\end{theorem}
\proof
\begin{enumerate}
	\item [(i)]
If $d=0$, then $G\cong K_n \cup K_1$ and if $d=n$, then $G\cong K_{n+1}$. So  the result for this case follows from Theorem \ref{hh}. Now let $1\leq d \leq n-1$  and denote the vertices of $K_n$ by $v_1, ... , v_n$ such that $w$ is adjacent to $v_1, ... , v_d$.\\
 If $1\leq d\leq n-1$, then $G$ has $d$ vertices of maximum degree $\Delta (G)=n$, namely the neighbors of $w$. Then, by removing  the vertex $w$, we have $G-w=K_n$ and maximum degree $K_n$ is equal with $n-1$. So $vs_{\Delta}(G)=1$.\\
For even  $n$,  we consider two cases.

\medskip
Case (a): If $1 \leq d \leq \frac{n}{2}$, then $G$ is in the class $1$. Consider the proper
edge coloring of $K_n$ with $n-1$ colors, where the vertices are in order 
$v_1, v_{d+1}, v_2,
v_{d+2}, . . . , v_{d-1},
 v_{2d-1},\\ v_{d}, v_{2d}$, $ v_{2d+1}, . . . , v_n$. 
Then the edges $v_1v_{d+1}, . . . , v_{d-1}v_{2d-1}$
are colored pairwise differently. Color these edges as well as edge $wv_d$ with the
new color $n$ and then color $wv_i$ with the old color of $v_iv_{d+i}$, $i = 1, ... , d-1$. This
implies $\chi^\prime (G)=\Delta (G)=n$.
With removing $w$, we have $K_n$ with $n$ even. Then $\chi^\prime (K_n)=\Delta (K_n)=n-1$. So $n=\chi^\prime (G) > \chi^\prime (K_n -w)=\chi^\prime (K_n)=n-1$ and therefore $vs_{\chi^\prime}(G)=1$.

\medskip
Case (b): If $\frac{n}{2} < d \leq n-1$, then $G$ is overfull, because  $|E(G)|=\frac{n(n-1)}{2}+d > \frac{n(n-1)}{2}+\frac{n}{2}=\frac{n^2}{2}=\Delta(G) (|V(G)|-1)/2$. Then $\chi^\prime (G)=\Delta (G)+1=n+1$. So $n+1=\chi^\prime (G) > \chi^\prime (G-w)=\chi^\prime (K_n)=n$. Therefore $vs_{\chi^\prime}(G)=1$.

\item[(ii)] 
 Since $n\geq 2$  is odd and $K_n \subseteq G \subseteq K_{n+1}$, $\chi^\prime (K_n)=\chi^\prime (G)=\chi^\prime (K_{n+1})=n$. By Corollary \ref{bbb} and Theorem \ref{hh}, $vs_{\chi^\prime}(G) \geq vs_{\chi^\prime}(K_n)=1$. If we remove $w$, then $G-w=K_n$. Therefore $vs_{\chi^\prime}(G)>1$. If we remove a vertex $v$ of $K_n$, then by proof of part (i), we have $vs_{\chi^\prime}(G-v)=1$ and therefore $vs_{\chi^\prime}(G)=2$.\qed
\end{enumerate}

\begin{figure}
	\begin{center}
		\includegraphics[width=0.5\textwidth]{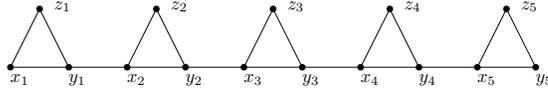}
		\caption{\label{fw} Graph $G_5$.}
	\end{center}
\end{figure}

\begin{theorem}
For every $k\in \mathbb{N}$, there exists a graph $G$ such that $vs_{\chi^\prime}(G)=k$. 
\end{theorem}
\proof
Let $G_k$, $k\geq 1$, be the graph of order $3k$ with $V(G_k )=\{x_1,y_1,z_1,...,x_k,y_k,z_k\}$, where the vertices $x_1, y_1,..., x_k, y_k$ (in that order) induce a path of length $2k-1$, and $z_i$ is adjacent to $x_i$ and $y_i$ for $1\leq i \leq k$ (see  $G_5$ in Figure \ref{fw}). Clearly $vs_{\chi^\prime}(G_k)=k$.\qed

\section{$\rho$-vertex stability for corona and join of two graphs} 
In this section,  we study the $vs_{\chi'}$ for corona product and join  of two graphs. First we consider the corona product of two graphs.  We recall that the corona of $G$ and $H$ is denoted by $G\circ H$, is a graph made by a copy of $G$ (which has $n$ vertices) and $n$ copy of $H$ and joining the $i$-th vertex of $G$ to every vertex in the $i$-th copy of $H$. Obviously  $\Delta (G\circ H)=\Delta(G)+|V(H)|$.

\begin{proposition}
If $\Delta (G)\geq 1$, then $vs_\Delta (G\circ H)\leq vs_\Delta (G)$.
\end{proposition}
\proof
Let $V^\prime$ be subset of vertices of $G$ which its removing  get $vs_\Delta(G)$. If  remove $V^\prime$ from  $G \circ H$, then $\Delta ((G \circ H)-V^\prime ) \leq \Delta (G-V^\prime )+n_2 < \Delta (G)+n_2=\Delta (G \circ H)$. Thus $vs_\Delta (G\circ H)\leq vs_\Delta (G)$.\qed
\begin{theorem}
Let $\Delta (G)\geq 1$ and $V^\prime$ be the subset of vertices of $G$ which its removing  get 
$vs_{\chi^\prime}(G)$. If $\Delta (G-V^\prime )+1 < \Delta (G)$, then $vs_{\chi^\prime}(G \circ H)=vs_{\chi^\prime}(G)$.
\end{theorem}
\proof
By Corollary \ref{bbb}, we have $vs_{\chi^\prime}(G) \leq vs_{\chi^\prime}(G\circ H)$. If remove $V^\prime$ from $G \circ H$, then $\chi^\prime ((G\circ H)-V^\prime) \leq \Delta ((G\circ H)-V^\prime)+1 \leq \Delta (G-V^\prime)+n_2+1 < \Delta (G)+n_2+1 \leq \chi^\prime (G \circ H)$. So $vs_{\chi^\prime} (G\circ H) \leq vs_{\chi^\prime} (G)$. Therefore $vs_{\chi^\prime} (G\circ H)=vs_{\chi^\prime} (G)$.\qed

\begin{theorem}\label{ii}
If $G\circ H$ is in the class $2$, then $vs_{\chi^\prime} (G\circ H)\leq \gamma (V_{\Delta (G)})$.
\end{theorem}
\proof
Since $G \circ H$ is in the class $2$, so $\chi^\prime (G\circ H)=\Delta (G \circ H)+1$. Let $\Gamma$ be a minimum dominating set of $V_{\Delta (G)}$. We remove the vertices of $\Gamma$ from  $G \circ H$. Since  $\Delta (G-\Gamma )<\Delta (G)$, so $\chi^\prime ((G \circ H)-\Gamma ) \leq \Delta ((G \circ H)-\Gamma )+1 \leq \Delta (G-\Gamma )+n_2+1 < \Delta (G)+n_2+1=\Delta (G \circ H)+1=\chi^\prime (G \circ H)$. Therefore $vs_{\chi^\prime} (G\circ H)\leq \gamma (V_{\Delta (G)})$.\qed

\begin{corollary}
Let $\Gamma$ be  a minimum dominating set of $V_{\Delta (G)}$. If $(G \circ H)-\Gamma$ is  in the class $1$, then $vs_{\chi^\prime} (G\circ H)\leq \gamma (V_{\Delta (G)})$.
\end{corollary}
\proof
Since $(G \circ H)-\Gamma$ is in the class $1$, so $\chi^\prime ((G \circ H)-\Gamma)=\Delta ((G\circ H)-\Gamma )$. Thus $\chi^\prime ((G \circ H)-\Gamma ) \leq \Delta ((G \circ H)-\Gamma ) \leq \Delta (G-\Gamma )+n_2 < \Delta (G)+n_2=\Delta (G \circ H)=\chi^\prime (G \circ H)$. Therefore the result is hold.\qed

\medskip
The join of two graphs $G_1$ and $G_2$, denoted by $G_1\vee G_2$, is a graph with vertex set $V(G_1)\cup V(G_2)$ and edge set $E(G_1)\cup E(G_2)\cup \{uv|u\in V(G_1) ~ and~ v\in V(G_2)\}$. 
Let  $n_1 = |V_1|$, $n_2 = |V_2|$, $\Delta _1=\Delta (G_1)$ and $\Delta _2 = \Delta (G_2)$. Clearly, $n(G_1\vee G_2) = n_1+n_2$ and $\Delta (G_1\vee G_2) = max \{n_1+\Delta _2, n_2+\Delta _1\}$.
\begin{proposition}
 $vs_\Delta (G_1\vee G_2) \leq 2$.
\end{proposition}
\proof
Let $v\in V(G_2)$ and  $u\in V(G_1)$. We have $\Delta (G-\{v,u\}) \leq max \{\Delta_1 +n_2-1, \Delta _2 +n_1-1\} < max \{\Delta_1 +n_2, \Delta _2 +n_1\}=\Delta (G)$. Therefore $vs_\Delta (G) \leq 2$.\qed
\begin{proposition}
If $G=G_1 \vee G_2$ and $\Delta _1+n_2 \neq \Delta _2+n_1$, then $vs_\Delta (G)=1$.
\end{proposition}
\proof
Without loss of generality, let $\Delta _1+n_2 > \Delta _2+n_1$. We remove a vertex $v$ of $G_2$. Then $\Delta (G-v) \leq \Delta (G_1)+n_2-1 < \Delta (G)$. So $vs_\Delta (G)=1$.\qed
\begin{theorem}
If $G=G_1 \vee G_2$ and $n_1,n_2\geq 2$, then $vs_{\chi ^\prime}(G) \leq 4$.
\end{theorem}
\proof
We consider the following cases:\\
Case 1) Let $\Delta _1+n_2=\Delta _2+n_1$. In this case, we remove a vertex $v$ of $G_2$ and a vertex $u$ of $G_1$. Then $\Delta (G-\{v,u\}) \leq max \{\Delta _1+n_2, \Delta _2+n_1\}$. Now, if $G-\{v,u\}$ is  in the class $1$, then $vs_{\chi ^\prime}(G) \leq 2$, and if $G-\{v,u\}$ is  in the class $2$, then we remove a vertex $v^\prime$ of $G_2$ in $G-\{v,u\}$ and a vertex $u^\prime$ of $G_1$ in $G-\{v,u\}$. Thus $\Delta (G-\{v,u,v^\prime ,u^\prime \}) \leq max \{\Delta _1+n_2-2, \Delta _2+n_1-2\}$. So $\chi ^\prime (G) < \chi ^\prime (G-\{v,u,v^\prime ,u^\prime \})$ and therefore in this case, $vs_{\chi ^\prime}(G) \leq 4$.\\
Case 2) Let $\Delta _1+n_2 \neq \Delta _2+n_1$. Without loss of generality, we  assume that $\Delta _1+n_2 > \Delta _2+n_1$.\\
Subcase a) If $G$ is in the class $2$, then $\chi ^\prime (G)=\Delta _1+n_2+1$. By removing a vertex $v$ of $G_2$,  $\Delta (G-v) \leq \Delta _1+n_2-1$. So $\chi ^\prime (G-v) \leq \Delta _1+n_2-1+1 < \Delta _1+n_2+1=\chi ^\prime (G)$. Therefore $vs_{\chi ^\prime}(G) =1$.\\
Subcase b) If $G$ is in the class $1$, then $\chi ^\prime (G)=\Delta _1+n_2$. We remove a vertex $v$ of $G_2$ and a vertex $u$ of $G_1$. So $\Delta (G-\{v,u\}) \leq max\{\Delta _1+n_2-1,\Delta _2+n_1-1\}$. If $G-\{v,u\}$ is  in class $1$, then $vs_{\chi ^\prime}(G) \leq 2$. Now let $G-\{v,u\}$ is in the class $2$. We remove a vertex $v^\prime$ of $G_2$ in $G-\{v,u\}$. Thus $\Delta (G-\{v,u,v^\prime \}) \leq max\{\Delta _1+n_2-2,\Delta _2+n_1-1\}$. So $\chi^\prime (G-\{v,u,v^\prime \}) \leq \Delta _1+n_2-1 < \Delta _1+n_2=\chi ^\prime (G)$. Therefore $vs_{\chi ^\prime}(G) \leq 3$.\qed
\begin{theorem}
If $G=G_1 \vee G_2$ is in class $2$, and $\Delta _1+n_2\neq \Delta _2+n_1$, then $vs_{\chi ^\prime}(G) =1$.
\end{theorem}
\proof
Since $G$ is in the class $2$, so $\chi^\prime (G)=\Delta _1+n_2+1$ if $\Delta _1+n_2>\Delta _2+n_1$ and $\chi^\prime (G)=\Delta _2+n_1+1$ if $\Delta _1+n_2<\Delta _2+n_1$. We assume that $\Delta _1+n_2>\Delta _2+n_1$. By removing  a vertex $v$ of $G_2$,  $\Delta (G-v) \leq \Delta _1+n_2-1$ and so $\chi ^\prime (G-v)\leq \Delta (G-v)+1\leq \Delta _1+n_2 < \Delta _1+n_2+1=\chi ^\prime (G)$. Therefore $vs_{\chi ^\prime}(G) =1$. If $\Delta _1+n_2<\Delta _2+n_1$, then the proof is similar to previous case.\qed

We need the following results to end the paper by presenting a sufficient condition for $vs_{\chi ^\prime}(G_1\vee G_2) =1$
\begin{theorem}{\rm \cite{e}}\label{t1}
	\begin{enumerate} 
		\item[(i)] 
		Let $G = G_1 \vee G_2$ be a join graph with $n_1 \leq n_2$. If $\Delta _1 > \Delta _2$, then $\chi^\prime (G)=\Delta (G)$.
		\item[(ii)] 
		Let $G = G_1 \vee G_2$ be a join graph with $n_1 = n_2$. If $\Delta _1 < \Delta _2$, then $\chi^\prime (G)=\Delta (G)$.
		\end{enumerate} 
\end{theorem}

\begin{theorem}
	\begin{enumerate} 
\item[(i)] If $G=G_1 \vee G_2$  with $n_1=n_2+1$ and $\Delta _1>\Delta _2$, then $vs_{\chi ^\prime}(G) =1$.
	\item[(ii)]  
	If $G=G_1 \vee G_2$ with $n_1=n_2+1$ and $\Delta _1<\Delta _2$, then $vs_{\chi ^\prime}(G) =1$.
\end{enumerate} 
\end{theorem}
\proof
	\begin{enumerate} 
		\item[(i)]
We remove a vertex $v$ of $G_2$. So $n_1\leq n_2-1$ and $\Delta _1>\Delta (G_2-v)$. By Theorem \ref{t1} (i), $G-v$ is in the class $1$, and so $\chi ^\prime (G-v)\leq \Delta _1+n_2-1<\Delta _1+n_2=\chi ^\prime (G)$. Therefore $vs_{\chi ^\prime}(G) =1$.
\item[(ii)] 
By removing  a vertex $u$ of $G_1$, in $G-u$, $|V(G-u)|=|V(G_2)|=n_2$ and $\Delta (G_1-u)<\Delta _2$. By Theorem \ref{t1} (ii), $G-u$ is in the class $1$ and $\chi ^\prime (G-u)=\Delta _2+n_2$. We know that $\chi ^\prime (G)>\Delta _2+n_2$. Therefore $vs_{\chi ^\prime}(G) =1$.\qed
\end{enumerate}

\end{document}